\documentclass[11pt]{amsart}

\usepackage{amsmath,amsfonts}
\usepackage{amsthm}
\usepackage{enumerate}

\newtheorem{theorem}{Theorem}[section]
\newtheorem{proposition}[theorem]{Proposition}
\newtheorem{lemma}[theorem]{Lemma}
\newtheorem{corollary}[theorem]{Corollary}

\theoremstyle{definition}
\newtheorem{definition}[theorem]{Definition}
\newtheorem{example}[theorem]{Example}

\theoremstyle{remark}
\newtheorem{remark}[theorem]{Remark}

\numberwithin{equation}{section}

\newcommand{\nC}{\mathbb C}
\newcommand{\nD}{\mathbb D}
\newcommand{\nM}{\mathbb M}

\newcommand{\cC}{{\mathcal C}}
\newcommand{\cD}{{\mathcal D}}
\newcommand{\cF}{{\mathcal D}}
\newcommand{\cL}{{\mathcal L}}
\newcommand{\cP}{{\mathcal P}}

\newcommand{\ds}{d_F}
\newcommand{\dm}{d_m}

\DeclareMathOperator{\diam}{diam} 
\DeclareMathOperator{\conv}{conv\!}

\begin{document}

\title[Perturbations of Roots of Polynomials]%
{Perturbations of Roots under \\ Linear Transformations of
Polynomials}

\author{Branko \'{C}urgus}
\address{Department of Mathematics, Western Washington
University, \newline \hspace*{3mm} Bellingham, WA 98225, USA}
\email{curgus@cc.wwu.edu}
\author{Vania Mascioni}
\address{Department of Mathematical Sciences, Ball State
University,\newline \hspace*{3mm} Muncie, IN 47306-0490, USA}
\email{vdm@cs.bsu.edu}


\date{}


\begin{abstract}
Let $\cP_n$ be the complex vector space of all polynomials of degree
at most $n$.  We give several characterizations of the linear
operators $T:\cP_n\rightarrow\cP_n$ for which there exists a
constant $C > 0$ such that for all nonconstant $f\in\cP_n$ there
exist a root $u$ of $f$ and a root $v$ of $Tf$ with $|u-v|\leq C$.
We prove that such perturbations leave the degree unchanged and, for
a suitable pairing of the roots of $f$ and $Tf$, the roots are never
displaced by more than a uniform constant independent on $f$.  We
show that such ``good'' operators $T$ are exactly the invertible
elements of the commutative algebra generated by the differentiation
operator. We provide upper bounds in terms of $T$ for the relevant
constants.
\end{abstract}

\maketitle

\section{Introduction} \label{intro}

Let $n$ be a positive integer, and denote by $\cP_n$ the
$(n+1)$-dimensional complex vector space of all polynomials of
degree at most $n$.  Let $T$ be a linear operator from $\cP_n$ to
$\cP_n$. In \cite{CM1} we proved that for each non-constant
polynomial $f \in \cP_n$ the polynomials $f$ and $\,Tf$ have at
least one common root if and only if $T$ is a non-zero constant
multiple of the identity on $\cP_n$.  In other words, if $T$ is not
a multiple of the identity, then there exists a polynomial $f \in
\cP_n$ such that $f$ and $\,Tf$ do not share any roots. A natural
question to ask is: How far apart are the roots of $\,Tf$ from the
roots of $f$?

This requires that we introduce a measure of distance between finite
subsets of the complex plane $\nC$.  In Section~\ref{sd} we
introduce four such distances, among which are two common ones:
$d_H,$ the Hausdorff distance and $d_F$, the Fr\'{e}chet distance.
The main result of this article is the characterization of the set
of those linear operators $T:\cP_n \to \cP_n$ for which there exists
a constant $C > 0$ such that for all $f\in\cP_n$, the distance
between the roots of polynomials $f$ and $\,Tf$ is at most $C$. Here
the distance can be any of the four distances that we introduce,
which implies that this set of ``good'' operators will turn out not
to depend on the distance used.

A simple example of a ``bad'' operator is the operator $R:\cP_n \to
\cP_n$ which changes the sign of the independent variable, defined
by
 \[
(Rf)(z) : = f(-z), \ \ \ z \in \nC, \quad f\in\cP_n.
 \]
If the roots of $f$ are large positive numbers, then the roots of
$Rf$ are negative numbers with large moduli, making any of the
distances that we consider as large as we want.

To illustrate our result let $\alpha$ be a complex number and
consider the linear operator $S(\alpha):\cP_n \to \cP_n$
corresponding to the additive shift of the independent variable. It
is defined by
\begin{equation} \label{eqSa}
\bigl(S(\alpha)f\bigr)(z) : = f(z+\alpha), \ \ \ z \in \nC,\quad
f\in\cP_n.
\end{equation}
It will be quite clear that (with respect to any of the four
distances) the distance between the roots of $f$ and the roots of
$S(\alpha)f$ will be at most $|\alpha|$ (see
Proposition~\ref{transl}).  The Taylor formula at $z$ implies that
the operator $S(\alpha)$ can be expressed as
 \[
S(\alpha) = I + \frac{\alpha}{1!}\, D +  \frac{\alpha^2}{2!}\, D^2 +
\cdots +  \frac{\alpha^n}{n!}\, D^n,
 \]
where $D:\cP_n \to \cP_n$ is the operator of differentiation with
respect to the complex variable.  This example hints at the main
result of this article stated in Theorems~\ref{tc1} and \ref{tc2}.
We paraphrase it below.

Let $T \in \cL(\cP_n), \, T\neq0$, and let $Z(f)$ denote the set of
the roots of a non-constant $f\in\cP_n$.  The following statements
are equivalent.
\begin{enumerate}[{\rm (i)}]
\item \label{sA}
There exists a constant $C > 0$, which depends on the distance $d$,
such that $d\bigl(Z(f),Z(Tf) \bigr) \leq C$ for each non-constant $f
\in \cP_n$.
\item \label{sB}
There exist $a_0, a_1, \ldots, a_n \in \nC, \, a_0 \neq 0$, such
that
\begin{equation} \label{eqinD}
T = a_0\, I + a_1\, D + a_2\, D^2 + \cdots + a_n\, D^n.
\end{equation}
\end{enumerate}
In (\ref{sA}) the symbol $d$ can be replaced with any of the
distances $d_m, d_h, d_H, d_F$ from Section~\ref{sd}. Thus
(\ref{sA}) really stands for four equivalent statements. Moreover,
for $T$ described in (\ref{sB}) and for each of the distances, in
Theorem~\ref{t13} we give an estimate for the maximum possible
distance between the roots of $f$ and the roots of $Tf$ in terms of
$T$.

Surprisingly, we found only one article, \cite{T},  which considers
the relationship between (\ref{sA}) and (\ref{sB}) as stated above.
In \cite{T} an entirely different method was used to prove that
(\ref{sB}) implies (\ref{sA}) with the distance $d = \ds$. The
converse was not considered in \cite{T}.  Also, no specific estimate
for $C$ is given there, which is in part due to the use of ``soft''
theorems from complex function theory.

It is not surprising, though, that the location of the roots of
$\,Tf$ in relation to the roots of $f$ for $T$ as in (\ref{sB}) has
been extensively researched, see \cite[Sections~5.3 and 5.4]{RS}. In
fact, the implication (\ref{sB})$\Rightarrow$(\ref{sA}), with the
distance $d = d_h$, is a consequence of Grace's theorem,
\cite[Theorem~5.3.1]{RS}.  For completeness we include the details
in Sections~\ref{sG} and \ref{sGT} below. Furthermore,
\cite[Corollary~5.4.1]{RS} is fundamental for the proof that
(\ref{sB}) implies (\ref{sA}) with the distance $d = \ds$.

The article is organized in twelve short sections; the first section
being this introduction.  In Section~\ref{sG} we recall Grace's
theorem and one of its consequences.  This consequence is of
interest to us since in Section~\ref{sGT} we restate it in terms of
operators on~$\cP_n$. To this end, in Section~\ref{scD} we study the
algebra of operators given by \eqref{eqinD} and a connection between
this algebra and $\cP_n$ is explored in Section~\ref{scDcP}. Since
the proofs in Sections~\ref{scD} and~\ref{scDcP} are short and
interesting we have not omitted them.  In Section~\ref{sGT} we
present a version of Grace's theorem for linear operators on
$\cP_n$. This theorem and a result from Section~\ref{sslo} are the
main tools in Section~\ref{smt1}, in which the first version of our
main result is formulated as Theorem~\ref{tc1}.  This is where we
prove the equivalence of (\ref{sA}) and (\ref{sB}) stated above. In
addition, Theorem~\ref{tc1}, among several equivalent statements,
contains the converse of Grace's theorem for linear operators. In
Theorem~\ref{tc1} we do not use the concepts of distances from
Section~\ref{sd}. We wanted to keep the first part of the article,
Sections~\ref{sG} through \ref{smt1}, independent of these concepts.

However, the second part depends heavily on the four distances from
Section~\ref{sd}. Why these four distances? The distance $\dm$ is
the simplest (``the two closest points distance"), $d_h$ and $d_F$
are implicitly already present in theorems about roots of
polynomials, and the Hausdorff distance $d_H$ is probably the
simplest  distance which is a metric. In Sections~\ref{sec} and
\ref{ses}, for each distance, we give exact calculations and
estimates for the maximum possible distance between the roots of $f$
and the roots of $\,Tf$ in terms of $T$. Finally, in
Section~\ref{smt2} we present the main theorem, Theorem~\ref{tc2}.
We conclude with several examples in Section~\ref{se}.

We now introduce the basic notation. By $\deg(f)$ we denote the
degree of a polynomial $f \in \cP_n$. For a non-zero polynomial $\,f
\in \cP_n$, $Z(f) \subset \nC$ will denote the multiset of all the
roots of $f$, that is, each root of $f$ appears in $Z(f)$ as many
times as its multiplicity as a root of $f$. Thus $Z(f)$ has exactly
$\deg(f)$ elements, and these are not necessarily distinct.  The
distinction between sets and multisets is essential only when we
consider the Fr\'{e}chet distance $d_F$.  In all other cases $Z(f)$
can be considered simply as the set of roots of $f$. For
completeness we set $Z(0) = \nC$.

By $\cL(\cP_n)$ we denote the set of all linear operators from
$\cP_n$ to $\cP_n$. We shall simply refer to elements in
$\cL(\cP_n)$ as operators.  Whenever we need a basis for $\cP_n$ we
shall use the basis $\bigl\{\phi_0,\phi_1,\ldots,\phi_n\bigr\}$ (in
this listed order) where $\phi_k(z) := z^k/k!, \,k = 0, 1, \ldots,
n$.

By $\nD(w,r)$ we denote the closed disk in $\nC$ centered at $w \in
\nC$ with radius $r > 0$. The letter $z$ always stands for a complex
number. For $A,B \subset \nC$ we define $A+B := \{u+v:u\in A,v\in
B\}$ and $-A:=\{-u:u\in A\}$.  By $\conv(A)$ we denote the convex
hull of $A$. Thus, for $f\in \cP_n$, $\conv\bigl(Z(f)\bigr)$ is the
convex hull of the roots of $f$.

Finally, we thought it would be in reader's interest to try to have
as many references as possible pointing to a single source, and the
recent monograph by Rahman and Schmeisser \cite{RS} is perfectly
adapted to the task. Another standard reference in this field is
\cite{M}.


\section{Grace's Theorem} \label{sG}

We begin with the definition of the \mbox{$*$-p}roduct of
polynomials which appears in \cite[page 375]{S} (see also  pages~148
and 178 in \cite{RS}).

\begin{definition} \label{defapol}
Let $f, g \in \cP_n$ be polynomials such that $\deg f = \deg g = m
> 0$. Let $a_k = f^{(k)}(0)$ and $b_k = g^{(k)}(0)$ for $k
=0,\ldots,n$, be the coordinates of $f$ and $g$ with respect to the
basis $\bigl\{\phi_0,\phi_1,\ldots,\phi_n\bigr\}$ of $\cP_n$. Set
\begin{equation*} 
(f * g)(z) := \sum_{k=0}^m b_k f^{(m-k)}(z) = \sum_{k=0}^m a_k
g^{(m-k)}(z) = \sum_{k=0}^{m} \Biggl(\!\sum_{j=0}^{m-k}
a_{k+j}\,b_{m-j}\!\Biggr)\frac{z^k}{k!} .
\end{equation*} 
\end{definition}
The reader can easily verify (or see \cite{S}) that the three sums
that appear in the definition are equal for all $z\in\nC$.

Definition~\ref{defapol} requires that the polynomials which are
being \mbox{$*$-mu}ltiplied have the same degree. The definition
depends on the common degree and the $*$-product is a polynomial of
the same degree.  It is clear that $*$ is commutative.

The next definition is equivalent to the standard one, see
\cite[Definition~3.3.1]{RS}. As before, $(Rf)(z) = f(-z)$.

\begin{definition}
Two polynomials $f$ and $g$ with equal positive degrees are {\em
apolar} if $\bigl(f * (Rg)\bigr)(0) = 0$.
\end{definition}

The symmetry of the apolarity relation follows from the
straightforward equality $g* (Rf) = (-1)^m f* (Rg)$.

The most important result about apolar polynomials is Grace's
theorem, see \cite[Theorem~3.4.1]{RS}.

\begin{theorem}[Grace]
Let $f$ and $g$ be apolar polynomials. If $\,\Omega$ is a circular
domain and $Z(g) \subset \Omega$, then $0 \in Z(f) - \Omega$.
\end{theorem}

With $S(\alpha)$ as defined in \eqref{eqSa} we clearly have
$Z\bigl(S(\alpha)p\bigr) = \{-\alpha\} +Z(p)$ for any polynomial
$p$. Let now $f$ and $g$ be polynomials with the same positive
degree. It is easy to verify that $ \bigl( S(\alpha)f\bigr)*(Rg)=
S(\alpha)\bigl(f*(Rg)\bigr). $ Combining the last two equalities we
conclude that $S(\alpha)f$ and $g$ are apolar if and only if $\alpha
\in Z(f*(Rg))$. Applying Grace's theorem to $S(\alpha)f$ and $g$
yields that $Z\bigl(f*(Rg)\bigr) \subset \ Z(f) - \Omega$ whenever
$\Omega$ is a circular domain and $Z(g) \subset \Omega$. Since
$Z(Rg) = -Z(g)$ this leads to the following theorem, see
\cite[Theorem~5.3.1]{RS}.
\begin{theorem} \label{tSt}
Let $f$ and $g$ be polynomials with the same positive degree. If
$\,\Omega$ is a circular domain and $Z(g) \subset \Omega$, then
$Z(f*g) \subset Z(f) + \Omega$.
\end{theorem}

For a fixed $g \in \cP_n$ with degree $n$, the mapping $f \mapsto f
* g$, $f\in\cP_n\!\setminus\!\cP_{n-1}$, is a restriction of a
linear combination of derivatives.  We shall explore the
relationship between the \mbox{$*$-p}roduct and operators on $\cP_n$
further. For that purpose we first study linear combinations of
derivatives.

\section{The commutative algebra $\cF(\cP_n)$} \label{scD}

The following definition introduces our main object of study.
\begin{definition} \label{fdef}
Let $D:\cP_n \to \cP_n$ be the operator of differentiation on
$\cP_n$. By $\cF(\cP_n)$ we denote the linear span in $\cL(\cP_n)$
of the operators $I, D, \ldots, D^n$. If $\alpha_0, \ldots, \alpha_n
\in \nC$ and
\begin{equation} \label{eqTc}
T = a_0\, I + a_1\, D +  \cdots + a_n\, D^n \in \cF(\cP_n),
\end{equation}
then we write $T = T(a_0, \ldots, a_n)$.
\end{definition}

To get familiar with the operators in $\cF(\cP_n)$ we first obtain
their matrix representation with respect to the basis of $\cP_n$
defined as
\begin{equation} \label{eqssB}
\bigl\{\phi_0,\phi_1,\ldots,\phi_n\bigr\} , \ \ \ \ \ \phi_k(z) :=
z^k/k!, \ \ \ k=0,\ldots,n.
\end{equation}
For each $m \in \{0,\ldots,n\}$ we clearly have
\begin{equation} \label{eqD}
D^k\phi_m = \phi_{m-k},  \ 0 \leq k \leq m, \ \ \ \text{and} \ \ \
D^k\phi_m = 0, \ m < k \leq n.
\end{equation}
Consequently, for $T$ given by \eqref{eqTc}, we have
\begin{equation} \label{eqD1}
T\phi_m = a_{m}\phi_{0} + a_{m-1}\phi_{1} + \cdots + a_1\phi_{m-1} +
a_0\phi_{m}, \ \ \ m=0,\ldots,n,
\end{equation}
and therefore $a_m = \bigl(T\phi_m\bigr)(0),\,m=0,\ldots,n$.
Equalities \eqref{eqD1} imply that the matrix of $\,T$ with respect
to the basis in \eqref{eqssB} of $\cP_n$ is the following upper
triangular Toeplitz matrix
 \begin{equation*} 
 \begin{bmatrix}
 a_0 & a_1 & a_2 & \cdots   & a_{n-1} & a_n \\[6pt] %
 0&    a_0 & a_1 & \cdots   & a_{n-2} & a_{n-1} \\[6pt] %
 0 &   0   & a_0 & \cdots   & a_{n-3} & a_{n-2}  \\[6pt] %
 \vdots & \vdots &  \ddots &
 \ddots &
   \ddots & \vdots \\[6pt]  %
 0 & 0 &  0 &  \cdots   & a_{0} & a_1 \\[6pt]  %
 0 & 0 & 0 & \cdots   & 0 & a_{0}
\end{bmatrix}.
\end{equation*}

Additional basic information about $\cF(\cP_n)$ is provided in the
next three statements.
\begin{proposition} \label{commutant}
Let $T \in \cL(\cP_n)$. Then $T\in \cF(\cP_n)$ if and only if $\,T$
commutes with $D$.
\end{proposition}
\begin{proof}
All elements of $\cF(\cP_n)$ clearly commute with $D$. To prove the
converse, set $\cC_n = \bigl\{T \in \cL(\cP_n)\, : \, TD=DT
\bigr\}$. Clearly $\cC_n$ is a subspace of $\cL(\cP_n)$ and
$\cF(\cP_n) \subset \cC_n$. Let $T \in \cC_n$. By \eqref{eqD}, $T
\phi_k = T D^{n-k}\phi_n = D^{n-k} T \phi_n$, for $k = 0,\ldots, n.$
Hence $T \in \cC_n$ is uniquely determined by $T\phi_n \in \cP_n$.
Consequently, the evaluation operator $T \mapsto T\phi_n, \,T \in
\cC_n$, is an injection. Therefore $\dim \cC_n \leq \dim \cP_n =
n+1$. Since $I,D,\ldots,D^n$, are linearly independent elements of
$\,\cC_n$, it follows that $\dim\bigl(\cC_n\bigr) = n+1$.
Consequently $\cC_n = \cF(\cP_n)$.
\end{proof}
\begin{corollary} \label{cca}
$\cF(\cP_n)$ is a maximal commutative subalgebra of $\cL(\cP_n)$.
\end{corollary}

The proposition below can be proved in different ways. See, for
example, the last paragraph in Section~\ref{scDcP}. We include this
proof since its method is also used in Section~\ref{ses}.

\begin{proposition} \label{invert}
Let $\,T \in \cF(\cP_n)$.  The operator $T$ is invertible if and
only if $\,T\phi_0 \neq 0$.  If $\,T$ is invertible, then $T^{-1}
\in \cF(\cP_n)$.
\end{proposition}
\begin{proof}
The ``only if" part of the first statement is obvious. To prove the
``if" part assume that $T\phi_0 \neq 0$. First note that since $D
\in \cF(\cP_n)$ is nilpotent, each operator $I-\gamma\,D \in
\cF(\cP_n)$ is invertible and
\begin{equation} \label{eqHgi}
(I - \gamma D)^{-1} = I + \gamma D + \cdots + \gamma^n D^n \in
\cF(\cP_n).
\end{equation}

Now let $T \in \cF(\cP_n)$ be given by \eqref{eqTc} and assume
$T\phi_0 \neq 0$. Then, by \eqref{eqD1}, $a_0 = T\phi_0 \neq 0$.
Following \cite[Section~5.4, p.~151]{RS}, let $\gamma_1, \ldots,
\gamma_n$ be the roots of
\begin{equation*} 
a_0\, z^n + a_1\, z^{n-1} + \cdots + a_{n-1}\,z + a_n
\end{equation*}
counted according to their multiplicities. Then clearly,
\begin{align*} 
a_0 + \cdots + a_n\, z^n & =
 a_0 z^n \prod_{j=1}^n\,\bigl(z^{-1}-\gamma_j\bigl)
 = a_0 \prod_{j=1}^n\,\left(1-\gamma_j z\right),
\end{align*}
and therefore,
\begin{equation} \label{ph21}
T = T(a_0,\ldots,a_n) = a_0 \prod_{j=1}^n\,
\bigl(I-\gamma_j\,D\bigr).
\end{equation}
As a product of invertible operators, $T$ is invertible. Since the
inverse of each of its invertible factors is in $\cF(\cP_n)$,
Corollary~\ref{cca} implies $T^{-1} \in \cF(\cP_n)$.
\end{proof}

\section{The algebra $\cF(\cP_n)$ and the vector space $\cP_n$ }
\label{scDcP}

In this section we explore the relationship between $\cF(\cP_n)$ and
$\cP_n$.

\begin{definition} \label{dvpi}
Define $\varpi: \cF(\cP_n) \rightarrow \cP_n$ by
\begin{equation*}
 \varpi(T) := T\phi_n, \ \ \ \  T \in \cF(\cP_n).
\end{equation*}
\end{definition}
If $T=T(a_0,\dots,a_n)$, then, by \eqref{eqD1} with $m=n$,
\begin{equation} \label{eqcvpi}
\varpi(T) = a_0\,\phi_n+a_1\,\phi_{n-1} + \cdots +a_{n-1}\,\phi_1
+a_n\,\phi_0.
\end{equation}
Results from Section~\ref{scD} and \eqref{eqcvpi} yield the
following proposition.
\begin{proposition} \label{pvpi}
The operator $\varpi$ is a linear bijection.  The image under
$\varpi$ of the set of all invertible operators in $\,\cF(\cP_n)$ is
the set of all polynomials of degree $n$.
\end{proposition}
Since $\cF(\cP_n)$ is a commutative algebra, it is natural to use
$\varpi$ to equip $\cP_n$ with an algebra structure. We do that
next.
\begin{definition} \label{dexs}
The  $\star$-product is defined on $\cP_n$ by
\begin{equation} \label{eqexs}
f \star g := \varpi\bigl(\varpi^{-1}(f) \varpi^{-1}(g) \bigr), \ \ \
\ \ f, g \in \cP_n.
\end{equation}
\end{definition}
The properties of $\cF(\cP_n)$ and $\varpi$ yield the following
corollary.
\begin{corollary} \label{cpst}
\begin{enumerate}[{\rm (a)}]
\item
The vector space $\cP_n$ equipped with the \mbox{$\star$-pr}oduct is
a commutative algebra with unit $\phi_n$.
\item
The operator $\varpi:\cF(\cP_n)\rightarrow\cP_n$ is an algebra
isomorphism.
\item
The $\star$-invertible polynomials are exactly the polynomials of
degree $n$.
\item
For each $m \in \{0,1,\ldots,n-1\}$, the set
\begin{equation} \label{eqsgm}
\Bigl\{f \in \cP_n: \deg f = n, \ f^{(k)}(0) = 0, \ k=0,\ldots,n-m-1
\Bigr\}
\end{equation}
is a $\star$-subgroup of the $\star$-group
$\cP_n\!\setminus\!\cP_{n-1}$.
\end{enumerate}
\end{corollary}
Now we are ready to establish the connection with the operators on
$\cP_n$.
\begin{proposition} \label{pTvp}
Let $T\in \cF(\cP_n)$. Then,
\begin{gather*}
Tf = f \star \bigl(\varpi(T)\bigr) = \bigl(\varpi(T)\bigr)\star f,
\ \ \ \ \ f\in\cP_n,  \\
(Tf)\star g = f\star (Tg) = T(f\star g), \ \ \ \ f,g\in\cP_n.
\end{gather*}
\end{proposition}
\begin{proof}
Let $T\in \cF(\cP_n)$ and $f \in \cP_n$. If $f = \varpi(V) =
V\phi_n$, then $\varpi^{-1}(f) = V$, and therefore
$\bigl(\varpi^{-1}(f)\bigr)(\phi_n) = f$. Using successively the
commutativity \mbox{of $\star$}, the fact that $\varpi$ is an
algebra isomorphism, the definition of $\varpi$, and the last
equality, we calculate
\begin{equation*}
f \star \bigl(\varpi(T)\bigr) =\varpi(T)\star f =  \varpi\bigl( T
\varpi^{-1}(f)\bigr) = T\bigl((\varpi^{-1}(f))(\phi_n)\bigr)
 =  T f.
\end{equation*}
Now the second claim follows from the associativity of \mbox{the
$\,\star$-p}roduct.
\end{proof}

The definition in \eqref{eqexs} is convenient since it emphasizes
the connection between $\cF(\cP_n)$ and $\cP_n$. However, the
formula for the \mbox{$\star$-pr}oduct in terms of the coordinates
with respect to $\{\phi_0,\ldots,\phi_n\}$ is also useful. Let $f,g
\in \cP_n$ and $a_k = f^{(k)}(0)$ and $b_k = g^{(k)}(0)$ for $k
=0,\ldots,n$. Now we first use \eqref{eqcvpi} to express
$\varpi^{-1}(f)$ and $\varpi^{-1}(g)$ in terms of the coordinates of
$f$ and $g$, then we calculate the composition
$\varpi^{-1}(f)\varpi^{-1}(g)$, and again, we use \eqref{eqcvpi} to
get
\begin{equation} \label{eqstp}
(f\star g)(z) =\sum_{k=0}^{n} \Biggl(\, \sum_{j=0}^{n-k}
a_{k+j}\,b_{n-j} \Biggr)\frac{z^k}{k!}.
\end{equation}
If $\deg f = \deg g = n$, a comparison of \eqref{eqstp} and
Definition~\ref{defapol} with $m = n$, yields that $f*g = f\star g$.
More generally, the following proposition holds.

\begin{proposition} \label{p*st}
Let $f,g \in \cP_n$, $\deg f = m$ and $\deg g = n$. Then
\begin{equation} \label{eq*st}
f*(D^{n-m}g)  = f \star g.
\end{equation}
For each $m \in \{1,\ldots,n\}$ the set
$\cP_m\!\setminus\!\cP_{m-1}$ with the \mbox{$*$-pr}oduct is a
commutative group.  The mapping $D^{n-m}$ restricted to the set
\eqref{eqsgm} is an isomorphism between the commutative group
\eqref{eqsgm} equipped with the \mbox{$\star$-pr}oduct and
$\cP_m\!\setminus\!\cP_{m-1}$ with the \mbox{$*$-product}.
\end{proposition}
\begin{proof}
Set $a_k = f^{(k)}(0)$ and $b_k = g^{(k)}(0)$ for $k =0,\ldots,n$.
Now regroup the terms in \eqref{eqstp} and use $a_{m+1} = \cdots =
a_n = 0$ to get a proof of \eqref{eq*st}:
\begin{equation*}
f\star g
 = \sum_{j=0}^{n} a_{j} g^{(n-j)} %
 = \sum_{j=0}^{m} a_{j} \bigl(D^{n-m}g\bigr)^{(m-j)} %
 = f* D^{n-m}g.
\end{equation*}

The mapping $D^{n-m}$ restricted to the set \eqref{eqsgm} is clearly
a bijection between that set and $\cP_m\!\setminus\!\cP_{m-1}$. Let
now $h$ and $g$ be polynomials in the set~\eqref{eqsgm}. Using
\eqref{eq*st} and Proposition~\ref{pTvp} we calculate
\begin{align*}
(D^{n-m}h)* (D^{n-m}g) = (D^{n-m}h) \star g = D^{n-m}(h \star g ).
\end{align*}
This proves the last claim of the proposition.
\end{proof}

The fact that $\cP_m\!\setminus\!\cP_{m-1}$ with the $*$-product is
a commutative group was proved in \cite{S}. Moreover, in \cite{S}
the reader can find a nice formula for the inverses.

\section{Grace's theorem for linear operators} \label{sGT}

The next theorem is a restatement of Theorem~\ref{tSt} in terms of
operators on $\cP_n$.  The role of $g$ in Theorem~\ref{tSt} is now
played by an invertible operator $T\in\cF(\cP_n)$. The union of the
sets $Z(T\phi_k),\, k=1,\ldots,n$, plays the role of $Z(g)$.
\begin{theorem} \label{tGT}
Let $T\in \cL(\cP_n)$. If $\,T\in\cF(\cP_n)$ and $\,T$ is
invertible, then $Z(Tf) \subset Z(f) + \Omega$ for all $f\in
\cP_n\!\setminus\!\{0\}$ and for all circular domains $\Omega$ such
that $Z(T\phi_k) \subset \Omega, \, k=1,\ldots,n$.
\end{theorem}
\begin{proof}
Let $\,T$ be an invertible operator in $\cF(\cP_n)$. For a constant
non-zero $f$ the theorem is obvious. Let $f$ be a non-constant
polynomial in $\cP_n$. Set $m = \deg f$. Then, by \eqref{eqD1} and
Proposition~\ref{invert}, $\deg(Tf) = m$ and
$\deg\bigl(\varpi(T)\bigr) = \deg\bigl(T\phi_n\bigr)=n$.
Propositions~\ref{pTvp} and \ref{p*st} yield
\begin{equation} \label{eqT*}
Tf = f\star \bigl(T\phi_n\bigr) = f*\bigl(D^{n-m} T\phi_n\bigr).
\end{equation}
Let $\Omega$ be a circular domain such that $Z(T\phi_k) \subset
\Omega$ for all $k=1,\ldots,n$. Since by Proposition~\ref{commutant}
and \eqref{eqD}, $D^{n-m} T\phi_n = T D^{n-m}\phi_n = T\phi_m$, the
theorem follows from \eqref{eqT*} and Theorem~\ref{tSt}.
\end{proof}

\begin{corollary} \label{cGT}
Let $\,T$ be an invertible operator in $\cF(\cP_n)$ and let $\Omega$
be a convex circular domain such that $Z(T\phi_n) \subset \Omega$.
Then $Z(Tf) \subset Z(f) + \Omega$ for all $f\in
\cP_n\!\setminus\!\{0\}$.
\end{corollary}
\begin{proof}
By Gauss-Lucas' theorem, see \cite[Theorem~2.1.1]{RS},
$Z\bigl(D^{n-k}T\phi_n\bigr) \subset
\conv\bigl(Z(T\phi_n)\bigr),\,k=1,\ldots,n$. Hence, for a convex
$\Omega$, $Z(T\phi_n) \subset \Omega$ implies $Z(T\phi_k)\subset
\Omega,\,k=1,\ldots,n$, and Theorem~\ref{tGT} applies.
\end{proof}

Theorem~\ref{tGT} is a motivation for the following definition.

\begin{definition} \label{dGS}
An operator $T\in\cL(\cP_n)$ will be called {\em Grace operator} if
there exists a finite set $A \subset \nC$ such that $Z(Tf) \subset
Z(f) + \Omega$ for all $f\in \cP_n\!\setminus\!\{0\}$ and for all
circular domains $\Omega$ such that $A \subset \Omega$.
\end{definition}

The question whether each Grace operator is an invertible operator
in $\cF(\cP_n)$ will be answered by Theorem~\ref{tc1}.

\section{The first step towards the main result} \label{sslo}

\begin{lemma} \label{l3}
Let $\,T \in \cL(\cP_n), \, T\neq 0$. Assume that there exists a
constant $C > 0$ such that for each non-constant polynomial $f \in
\cP_n$ there exist $u \in Z(f)$ and $v \in Z(Tf)$ such that $|u - v|
\leq C$. Then the matrix of $\,T$ with respect to the basis
$\bigl\{\phi_0,\ldots,\phi_n\bigr\}$ of $\,\cP_n$ is upper
triangular and the main diagonal entries are all equal to the same
non-zero constant $T \phi_0$.
\end{lemma}

\begin{remark} \label{rnc}
Note that the hypothesis of Lemma~\ref{l3} implies that $Z(Tf) \neq
\emptyset$ for all non-constant $f \in \cP_n$.  Further, the
conclusion of the lemma implies that $T$ maps constant polynomials
into constants.
\end{remark}

\begin{proof}[Proof of Lemma {\rm \ref{l3}}]
Let $m \in  \{ 1, \ldots, n\}$ and $t > 0$ be arbitrary. Consider
the polynomials $\phi_m - (t^m/m!) \phi_0$ and $T\bigl(\phi_m -
(t^m/m!)\phi_0\bigr)$. By hypothesis these two polynomials have
roots which are at most $C$ apart. Thus, for each $t > 0$ there
exists an $m$-th root of unity $\theta(t)$ such that the polynomial
$T \phi_m - (t^m/m!) T \phi_0$ has a root $w(t)$ in the disc
$\nD(t\, \theta(t),C)$.  If we assume that $T\phi_0 = 0$, then the
last statement would imply $T\phi_m = 0$. Since $m \in
\{1,\ldots,n\}$ is arbitrary, this would yield $T=0$. But $T\neq0$;
hence, $T\phi_0 \neq 0$ holds.

Set $v(t) = w(t)-t\,\theta(t), \, t > 0$. Since $v(t) \in \nD(0,C)$
and  $\theta(t) \in \nD(0,1)$, for each $k=0,\ldots,n$,
\begin{equation} \label{eqO}
(T \phi_k)\bigl(t\,\theta(t)+v(t)\bigr) =
O\bigl(t^{\deg(T\phi_k)}\bigr), \ \ \ t \to +\infty,
\end{equation}
and $\deg(T\phi_k)$ is the smallest power of $t$ for which
\eqref{eqO} holds. The special case of \eqref{eqO}, with  $k=0$,
implies, for each $m=0,\ldots,n$,
\begin{equation} \label{eqO0}
(t^m/m!) (T \phi_0)\bigl(t\,\theta(t)+v(t)\bigr) =
O\bigl(t^{m+\deg(T\phi_0)}\bigr), \ \ \ t \to +\infty,
\end{equation}
and $m+\deg(T\phi_0)$ is the smallest power of $t$ for which
\eqref{eqO0} holds. Recall that by the definition of $\theta(t)$ and
$v(t)$ we have
\begin{equation} \label{eqm0}
(T \phi_m)\bigl(t\,\theta(t)+v(t)\bigr) = (t^m/m!) (T
\phi_0)\bigl(t\,\theta(t)+v(t)\bigr), \ \ \ t > 0.
\end{equation}
This, \eqref{eqO} and \eqref{eqO0} imply
\begin{equation*}
\deg(T\phi_m) = m + \deg(T\phi_0), \ \ \ m = 0,1,\ldots,n.
\end{equation*}
Since $\deg(T\phi_n) \leq n$, the last equality with $m=n$  implies
$\deg(T\phi_0) = 0$; that is, $T\phi_0$ is constant. Consequently,
\begin{equation} \label{eqAem}
\deg(T\phi_m) = m, \ \ \ m = 0,1,\ldots,n.
\end{equation}
Hence, the matrix of $\,T$ with respect to the basis
$\bigl\{\phi_0,\ldots,\phi_n\bigr\}$ of $\,\cP_n$ is upper
triangular.

The main diagonal entries of this matrix are equal to the limits
 \[
\lim\limits_{z\to \infty} \frac{(T \phi_m)(z)}{\phi_m(z)}, \ \ \
m=0,1,\ldots,n,
 \]
existence of which is a consequence of \eqref{eqAem}. Since $v(t)
\in \nD(0,C)$ and  $\theta(t) \in \nD(0,1)$, \eqref{eqm0} implies
that these limits all equal to the constant $T\phi_0$.  The lemma is
proved.
\end{proof}
The following theorem is the first step towards a complete answer to
the question posed in the Introduction.
\begin{theorem} \label{tpa1}
Let $\,T \in \cL(\cP_n),\, T\neq 0$. Assume that there exists a
constant $C > 0$ such that for each non-constant $f \in \cP_n$ there
exist $u \in Z(f)$ and $v \in Z(Tf)$ such that $|u - v| \leq C$.
Then
\begin{equation} \label{eqsfA}
T = a_0\, I + a_1\, D + \cdots + a_n\, D^n,
\end{equation}
where
\begin{equation} \label{eqda}
a_0 = (T \phi_0)(0) \neq 0 \ \ \ \text{and} \ \ \ a_k = (T
\phi_k)(0), \ \ \ k = 1,\ldots,n.
\end{equation}
\end{theorem}
\begin{proof}
Let $T$ be as in the hypothesis. Define the coefficients $a_k, \, k
= 0,\ldots,n$, by \eqref{eqda}. Next we shall prove \eqref{eqsfA}.

From Lemma~\ref{l3} we know that the matrix of $T$ with respect to
the basis $\bigl\{\phi_0,\ldots,\phi_n\bigr\}$ of $\cP_n$ is upper
triangular and the main diagonal entries are all equal to the same
non-zero number $a_0 :=(T \phi_0)(0)$ (and thus the statement $a_0
\neq 0$ in \eqref{eqda} is justified).  This implies that $T$ is
invertible and
\begin{equation} \label{eqm01}
T\phi_0 =\bigl((T \phi_0)(0)\bigr) \phi_0 = a_0 \phi_0.
\end{equation}
(Remember that $\phi_0(z) := 1$ for all $z \in \nC$).

Since $\bigl\{\phi_0,\ldots,\phi_n\bigr\}$ is a basis for $\cP_n$,
\eqref{eqsfA} is equivalent to
\begin{equation} \label{AincF2}
T\phi_m = a_0\, \phi_m + a_1\, \phi_{m-1} + \cdots +  a_{m-1}\,
\phi_{1} + a_{m}\, \phi_0, \ \ \ m=0,\ldots,n.
\end{equation}

We prove \eqref{AincF2} by complete induction with respect to $m$.
By \eqref{eqm01} equality \eqref{AincF2} holds for $m = 0$.  Let
$k\in \{1,\dots, n\}$ and assume that \eqref{AincF2} is true for $m
= 0, \ldots, k-1$. We need to prove that \eqref{AincF2} holds for
$m=k$. By Lemma~\ref{l3}, $T \phi_k$ is a polynomial of degree $k$:
\begin{equation*}
\bigl(T \phi_k\bigr)(z) = \sum_{j=0}^k b_{k-j}\,\phi_j(z)
\end{equation*}
where $b_j \in \nC, \, j = 0,\ldots,k$. The rest of the proof is
devoted to calculating the coefficients $b_j, \, j=0,\ldots,k$.

Let $w \in \nC$ be arbitrary and consider the polynomial $p_k(z) :=
\phi_k(z - w)$.  Using the binomial expansion of $(z - w)^k$ and the
induction hypothesis, for all $z \in \nC$, we get the identity
 \allowdisplaybreaks{%
\begin{align*} 
\bigl(T p_k\bigr)(z) &= 
\sum_{j=0}^k\,\frac{1}{j!}(-w)^{j} \bigl(T \phi_{k-j}\bigr)(z) \\
 &= 
 \bigl(T \phi_k\bigr)(z) + \sum_{j=1}^{k}\,\frac{1}{j!}(-w)^{j}
    \sum_{l=0}^{k-j}\,a_{l} \phi_{k-j-l}(z) \\
 &= 
b_k + \sum_{l=0}^{k-1} b_{l} \, \phi_{k-l} +
\sum_{j=1}^{k}\,\frac{1}{j!}(-w)^{j}
    \sum_{l=0}^{k-j}\,\frac{a_{l}}{(k-j-l)!} z^{k-j-l} \\
 &= 
b_k + \sum_{l=0}^{k-1} \frac{b_{l}}{(k-l)!} z^{k-l}\\
& \phantom{= b_k \ } + \sum_{l=0}^{k-1} \frac{a_{l}}{(k-l)!} \left(
\sum_{j=1}^{k-l}\,
\binom{k-l}{j} \left(-\frac{w}{z}\right)^{j} z^{k-l} \right) \\
  &= 
b_k + \sum_{l=0}^{k-1}\, \left[\frac{b_{l}}{(k-l)!} +
\frac{a_{l}}{(k-l)!} \, \sum_{j=1}^{k-l}\,
 \binom{k-l}{j}
 \left(-\frac{w}{z}\right)^{j} \right]  z^{k-l}.
\end{align*}}
By hypothesis, for each $w \in \nC$ there exists $u(w) \in Z(T p_k)$
such that \mbox{$|w-u(w)|\leq C$}. Put $v(w) = u(w) - w$ and note
that $|v(w)| \leq C$ for all $w \in \nC$.  The substitution $z =
u(w)=w + v(w)$ in the last long displayed identity yields
\begin{equation*} 
b_k + \sum_{l=0}^{k-1}\frac{1}{(k-l)!}\!\left[b_{l}+ a_{l}
\sum_{j=1}^{k-l}\,\binom{k-l}{j}\!
 \left(\!\frac{-w}{w+v(w)}\!\right)^j \right]\!
 \bigl(w+v(w)\bigr)^{k-l} = 0,
\end{equation*}
which simplifies to
\begin{equation} \label{long2}
b_k + \sum_{l=0}^{k-1} \frac{1}{(k-l)!}\!\left[b_{l}+ a_{l}\!
  \left(\!\!\left(\frac{v(w)}{w+v(w)}\right)^{\!k-l}-1\!\right)\!\right]\!
 \bigl(w+v(w)\bigr)^{\!k-l} = 0,
\end{equation}
using
 \[
\sum_{j=1}^{k-l}\,\binom{k-l}{j}
 \left(- \frac{w}{w+v(w)}\right)^{\!j} = \left(1-
 \frac{w}{w+v(w)}\right)^{\!k-l} - 1.
 \]
Regrouping terms in \eqref{long2} yields
\begin{equation} \label{long3}
 b_k  +  \sum_{l=0}^{k-1} \frac{b_{l} - a_{l}}{(k-l)!}
   \,  \bigl(w+v(w)\bigr)^{k-l}
 + \sum_{l=0}^{k-1}  \frac{a_{l}}{(k-l)!}
 v(w)^{k-l}
  = 0.
\end{equation}
Since $|v(w)| \leq C$ for all $w \in \nC$ the last sum in
\eqref{long3} is a bounded function of $w$.  Therefore \eqref{long3}
implies
\begin{equation*} 
 \sum_{l=0}^{k-1} \, \frac{b_{l} - a_{l}}{(k-l)!}
   \,  \bigl(w+v(w)\bigr)^{k-l} = O(1), \ \ \ \
    |w| \to +\infty.
\end{equation*}
Again, since $|v(w)| \leq C$ for all $w \in \nC$, the last displayed
relation yields
\[
b_{l} - a_{l}= 0, \ \ \ l = 0, 1, \ldots, k-1.
\]
Since clearly $b_k = (T\phi_k)(0) = a_k$, we have proved that
\eqref{AincF2} holds for $m = k$. By induction, \eqref{AincF2} holds
for all $m = 0, 1, \ldots, n$, and the theorem is proved.
\end{proof}

\section{The first version of the main theorem} \label{smt1}

Let $f$ be a non-constant polynomial. The number
\[
\varrho[f] := \max\bigl\{|u|:u \in Z(f)\bigr\}
\]
is called the {\em root radius} of $f$.

\begin{proposition} \label{ppa2}
Let $\,T$ be an in invertible operator in $\cD(\cP_n)$. Then for
every non-constant $f \in \cP_n$ and for each $v \in Z(Tf)$ there
exists $u \in Z(f)$ such that $|v-u| \leq \varrho[T\phi_n]$.
\end{proposition}
\begin{proof}
Since the proposition is trivial for a non-zero constant multiple of
the identity operator, we assume that $T \neq a_0 I$. Then, by
\eqref{eqD1}, $\varrho[T\phi_n] > 0$.  Let $f\in \cP_n$ be a
non-constant polynomial. Since $\nD\bigl(0,\varrho[T\phi_n]\bigr)$
is a convex circular domain and $Z(T\phi_n) \subset
\nD\bigl(0,\varrho[T\phi_n]\bigr)$, Corollary~\ref{cGT} yields
\begin{equation*} 
Z(Tf) \subset  Z(f) + \nD\bigl(0,\varrho[T\phi_n]\bigr).
\end{equation*}
Thus, for each non-constant $f \in \cP_n$, for every $v \in Z(Tf)$
there exists $u \in Z(f)$ such that $|u-v| \leq \varrho[T\phi_n]$.
\end{proof}

\begin{remark} \label{rgC}
The conclusion of Proposition~\ref{ppa2} can also be expressed as
\[
Z(Tf) \subset \bigcup_{u\in Z(f)} \nD\bigl(u,\varrho[T\phi_n]\bigr).
\]
\end{remark}

The following theorem is the main result of this section.

\begin{theorem} \label{tc1}
Let $\,T \in \cL(\cP_n), \, T\neq 0$. The following statements are
equivalent.
\begin{enumerate}[{\rm (a)}]
\item \label{tc1c}
$T$ is an invertible operator in $\,\cF(\cP_n)$.
\item \label{tc1b}
$T$ is a Grace operator.
\item \label{tc1d}
There exists a constant $C_1>0$ such that for each non-constant $f
\in \cP_n$ and for each $w \in Z(Tf)$ there exists $v \in Z(f)$ such
that $|w-v| \leq C_1$.
\item \label{tc1e}
There exists a constant $C_2 > 0$ such that for each non-constant $f
\in \cP_n$ there exist $u \in Z(f)$ and $v \in Z(Tf)$ such that $|v
- u| \leq C_2$.
\item \label{tc1f}
$T$ is invertible and  it commutes with the differentiation operator
$D$.
\end{enumerate}
\end{theorem}
\begin{proof}
With $A = \bigcup \bigl\{Z(T\phi_k): k =1,\ldots,n\bigr\}$, the
implication (\ref{tc1c})$\Rightarrow$(\ref{tc1b}) follows from
Theorem~\ref{tGT}.

The following short proof of (\ref{tc1b})$\Rightarrow$(\ref{tc1d})
is similar to the proof of Proposition~\ref{ppa2}.  Assume
(\ref{tc1b}) and let $A \subset \nC$ be a finite set from
Definition~\ref{dGS}. Let $C_1 > 0$ be such that $A \subset
\nD(0,C_1)$.  Let $f\in \cP_n\!\setminus\!\cP_0$.  Since
$\nD(0,C_1)$ is a circular domain, by Definition~\ref{dGS} we  have
$Z(Tf) \subset Z(f) + \nD(0,C_1)$.  Hence,  for each $v \in Z(Tf)$
there exists $u \in Z(f)$ such that $v-u \in \nD(0,C_1)$. This
proves (\ref{tc1d}).

The implication (\ref{tc1d})$\Rightarrow$(\ref{tc1e}) is obvious and
(\ref{tc1e})$\Rightarrow$(\ref{tc1c}) was proved in
Theorem~\ref{tpa1}. Since (\ref{tc1c})$\Leftrightarrow$(\ref{tc1f})
is an immediate consequence of Proposition~\ref{commutant}, the
theorem is proved.
\end{proof}

\section{Distances} \label{sd}

Recall the following standard definition:
\begin{definition} \label{dm}
A function $d: X \times X \to [0,+\infty)$ is a {\em metric} on a
nonempty set $X$, if for all $x,y,z \in X$ we have
\begin{enumerate}[(a)]
\item \label{ida}
$d(x,x) = 0$;
\item \label{idb}
$d(x,y) = 0$ implies that $x = y$;
\item \label{idc}
$d(x,y) = d(y,x)$;
\item \label{idd}
$d(x,z) \leq d(x,y)+d(y,z)$.
\end{enumerate}
\end{definition}
The problem of measuring the distance between two finite sets of
points has been considered in several seemingly unrelated areas of
research. For a recent account and references see \cite{EM}.

The best known metric on the family of finite nonempty subsets of
$\nC$ is the {\em Hausdorff metric} defined as
\begin{equation*}
d_H(A,B) := \max\Bigl\{\max_{x\in A} \min_{y\in B} |x-y|, \
\max_{x\in B} \min_{y\in A} |x-y|\Bigr\},
\end{equation*}
where $A$ and $B$ are nonempty finite subsets of $\nC$. That the
function $d_H$ is really a metric on the family of finite nonempty
subsets of $\nC$ is a simple exercise.

To connect the Hausdorff metric to Theorem~\ref{tc1} we introduce
two related functions. We informally call these functions {\em
distances} since they indicate the location of points of one set in
relation to the other.  For two nonempty finite subsets $A$ and $B$
of $\nC$ define
\begin{align*}
d_m(A,B) &:= \min\bigl\{  |x-y| \; : \; x \in A, \; y\in B\bigr\},
\\
d_h(A,B) &:= \max\bigl\{ d_m\bigl(A,\{y\}\bigr) \; : \; y \in B
\bigr\}.
\end{align*}
Now, the Hausdorff metric can be expressed as
 \begin{equation*} 
d_H(A,B) = \max\bigl\{d_h(A,B),d_h(B,A) \bigr\}.
 \end{equation*}
Clearly neither of the functions $d_m$ and $d_h$ is a metric. The
function $d_m$ satisfies only (\ref{ida}) and (\ref{idc}) and $d_h$
satisfies only (\ref{ida}) and (\ref{idd}) in Definition~\ref{dm}.
The function $d_h$ is sometimes called asymmetric Hausdorff
distance. Both distances $d_m$ and $d_h$ appear implicitly in
Theorem~\ref{tc1}.

The definitions of $d_m$ and $d_h$ can be extended to include the
empty set and the entire complex plane, which correspond to $Z(f)$
for $f \in \cP_0\!\setminus\!\{0\}$ and $Z(0)$. For both $d=d_m$ and
$d=d_h$, we set $d(\emptyset,\emptyset) = 0$ and $d(A,\emptyset) =
d(\emptyset,A) = +\infty$ whenever $A\neq\emptyset$.  For either
$A=\nC$ or $B=\nC$ the original definitions make sense, giving the
value $0$. It is in this extended sense that $d_m, d_h$ and $d_H$
will be used in the rest of the article.

Another well known metric is the Fr\'{e}chet metric (see \cite{F}
where a similar definition was first introduced, or \cite[Chapter
6]{E} where an analogous metric is defined for curves). Let $m$ be a
positive integer and put $\nM = \{1,\ldots,m\}$.  By $\Pi_m$ we
denote the set of all permutations of $\nM$.  For two functions $u,
v :\nM \to \nC$ we define
\begin{equation*}
d_F(u,v) := \min_{\sigma \in \Pi_m} \max_{k\in\nM}\, \bigl|u(k) -
v(\sigma(k)) \bigr|.
\end{equation*}
The function $d_F$ is not a metric on $\nC^m$ since it does not
satisfy (\ref{idb}) in Definition~\ref{dm}.  But $d_F$ is a metric
on the factor set $\nC^m\!/\!\!\sim$, where $u \sim v\Leftrightarrow
d_F(u,v) = 0$.  The elements of the factor set $\nC^m\!/\!\!\sim$
can be identified with ``unordered'' $m$-tuples, that is with
multisets of $\,m$ complex numbers in which a same element can
appear more than once.  This concept fits well with the sets of
roots of polynomials where roots can occur with multiplicities.  In
this context the Fr\'{e}chet metric is defined for two multisets of
$m$ complex numbers $U = \{u_1,\dotsc,u_m\}$ and $V =
\{v_1,\dotsc,v_m\}$ by
\begin{equation*}
d_F(U,V) := \min_{\sigma \in \Pi_m} \max_{k\in\nM}\, \bigl|u_k -
v_{\sigma(k)} \bigr|.
\end{equation*}
In some ways this distance is a natural distance when perturbation
of the roots of polynomials are studied, see \cite[Theorem,
p.~276]{O}, \cite{K} and \cite{CM2}.  For a simple proof that $d_F$
is a metric see \cite{CM2}.

We are interested in the question of what happens (in a quantitative
sense) to the roots of polynomials under linear operators on
$\cP_n$. The following numbers give a ``one number summary'' answer
to this question for $T \in \cL(\cP_n)$:
\begin{align*}
K_m(T) &:= \sup\bigl\{d_m\bigl(Z(f),Z(Tf)\bigr): f \in
\cP_n\bigr\} ,  \\
K_h(T) &:= \sup\bigl\{d_h\bigl(Z(f),Z(Tf)\bigr): f \in
\cP_n\bigr\} ,   \\
K_H(T) &:= \sup\bigl\{d_H\bigl(Z(f),Z(Tf)\bigr): f \in \cP_n\bigr\}.
\end{align*}

The power of these definitions is in the fact that the statements
from theorems in Sections~\ref{sslo} and \ref{smt1} can now be
formulated in a more compact way.  For example the hypothesis of
Theorem~\ref{tpa1} is: $K_m(T) < +\infty$ and the conclusion of
Proposition~\ref{ppa2} is: $K_h(T) \leq \varrho[T\phi_n]$.

Since the Fr\'{e}chet metric is defined only for multisets with the
same number of elements, we can define $K_F(T)$ only for $T \in
\cL(\cP_n)$ with the property $\deg(Tf) = \deg(f)$ for every $f \in
\cP_n$. For such $T$ we define
\begin{equation*}
K_F(T) := \sup \bigl\{d_F\bigl(Z(f),Z(Tf)\bigr): f \in
\cP_n\!\setminus\!\cP_0 \bigr\}.
\end{equation*}

\begin{proposition} \label{ps}
Let $T, T_1 \in \cL(\cP_n)$ and let $\alpha$ be a non-zero complex
number. Then
\begin{enumerate}[{\rm (a)}]
\item \label{psc}
$K_m(T) \leq K_h(T) \leq K_H(T)$.
\item \label{psa}
$K_h(T_1 T) \leq K_h(T) + K_h(T_1)$.
\item \label{psen}
$K_H(T_1 T) \leq K_H(T) + K_H(T_1)$.
\item \label{psfn}
For invertible $T$, $K_H(T)  = \max\bigl\{K_h(T), K_h(T^{-1})\bigr\}
= K_H(T^{-1})$.
\end{enumerate}
If $\,\deg(f) = \deg(Tf) = \deg(T_1f)$ for all $f \in \cP_n$, then
\begin{enumerate}[{\rm (a)}]
 \setcounter{enumi}{4}
\item \label{pse}
$K_F(T_1 T) \leq K_F(T) + K_F(T_1)$.
\item \label{psg}
$K_H(T) \leq K_F(T)$.
\item \label{psh}
$\,T$ is invertible and $K_F(T) = K_F(T^{-1})$.
\end{enumerate}

\end{proposition}
\begin{proof}
The proofs follow directly from the definitions, properties of the
underlying distances and properties of the supremum. A proof of
(\ref{psa}) follows:
\begin{align*}
K_h(T_1 T) &= \sup\bigl\{d_h\bigl(Z(f),Z(T_1 Tf)\bigr): f \in
\cP_n\} \\
 & \leq \sup\bigl\{d_h\bigl(Z(f),Z(Tf)\bigr): f \in
\cP_n\bigr\} \\
& \phantom{\leq K_h(T)} \;
     + \sup\bigl\{d_h\bigl(Z(Tf),Z(T_1 Tf)\bigr):f \in
\cP_n\bigr\} \\
& \leq K_h(T) + K_h(T_1).
\end{align*}
To prove (\ref{psfn}) assume that $T$ is invertible. Then
$T^{-1}\cP_n = \cP_n$ and therefore
\begin{align*} 
K_h(T^{-1}) & = \sup\bigl\{d_h\bigl(Z(f),Z(T^{-1}f)\bigr): f \in
\cP_n\bigr\} \\
&= \sup\bigl\{d_h\bigl(Z(Tg),Z(g)\bigr): g \in \cP_n\bigr\}.
\end{align*}
Now the first equality in (\ref{psfn}) follows from the definition
of $K_H(T)$. The second equality follows from the first when $T$ is
substituted by $T^{-1}$.

The remaining statements are proved similarly.
\end{proof}

\section{Exact Calculations} \label{sec}

We start with general results for $K_h$ and $K_H$.

\begin{theorem} \label{tpa2K}
Let $\,T$ be an invertible operator in $\cF(\cP_n)$. Then
\begin{equation*} 
K_h(T) = \varrho[T\phi_n] = d_h\bigl(Z(\phi_n),Z(T\phi_n)\bigr).
\end{equation*}
\end{theorem}
\begin{proof}
By Proposition~\ref{ppa2},
\[
d_h\bigl(Z(f),Z(Tf)\bigr) \leq \varrho[T\phi_n], \ \ \ \ f \in
\cP_n\!\setminus\!\cP_0.
\]
Since $T$ maps constants onto constants, it follows that $K_h(T)
\leq \varrho[T\phi_n]$. Clearly,
$d_h\bigl(Z(\phi_n),Z(T\phi_n)\bigr) = \varrho[T\phi_n]$, and
therefore, $K_h(T) \geq \varrho[T\phi_n]$.
\end{proof}
Proposition~\ref{ps}(\ref{psfn}) now yields:
\begin{corollary} \label{cKHe}
Let $\,T$ be an invertible operator in $\cF(\cP_n)$. Then
\begin{align*}
 K_H(T) & = \max\bigl\{\varrho[T\phi_n],
 \varrho[T^{-1}\phi_n]\bigr\} \\
   & = \max\bigl\{
d_H\bigl(Z(\phi_n),Z(T\phi_n)\bigr),
d_H\bigl(Z(T^{-1}\phi_n),Z(\phi_n)\bigr) \bigr\}.
\end{align*}
\end{corollary}

\begin{remark}
Theorem~\ref{tpa2K} conveys that the worst possible perturbation of
roots measured by the distances $d_h$ occurs at the polynomial
$\phi_n$.  That is, for all $f \in \cP_n$,
\[
d_h\bigl(Z(f),Z(Tf)\bigr) \leq d_h\bigl(Z(\phi_n),Z(T\phi_n)\bigr) .
\]
For the distance $d_H$, by Corollary~\ref{cKHe}, the worst possible
perturbation of roots occurs either at $\phi_n$ or at
$T^{-1}\phi_n$. That is, for all $f \in \cP_n$,
\[
d_H\bigl(Z(f),Z(Tf)\bigr)\leq \max\bigl\{
d_H\bigl(Z(\phi_n),Z(T\phi_n)\bigr),
d_H\bigl(Z(T^{-1}\phi_n),Z(\phi_n)\bigr) \bigr\}.
\]
It would be interesting to know whether the last two inequalities
must be strict when $f$ is not a multiple of $\phi_n$.
\end{remark}

Exact calculations are possible only for simple operators in
$\cF(\cP_n)$. We study two such classes. Let $\alpha, \gamma
\in\nC$. As before, $S(\alpha) \in \cF(\cP_n)$ is the operator that
shifts the independent variable by $\alpha$ defined in \eqref{eqSa}.
Further, we define
\begin{equation*}
H_k(\gamma) := I -\gamma D^k \, ,  \ \ \ k = 1,\dotsc,n .
\end{equation*}

\begin{proposition} \label{transl}
Let $\alpha\in\nC$ and consider $S(\alpha) \in \cF(\cP_n)$. Then
\[
K_m\bigl(S(\alpha)\bigr) =  K_h\bigl(S(\alpha)\bigr) =
K_H\bigl(S(\alpha)\bigr) = K_F\bigl(S(\alpha)\bigr) = |\alpha|.
\]
\end{proposition}
\begin{proof}
For $f\in\cP_n$ we have $Z\bigl(S(\alpha)f\bigr) =\{-\alpha\}+Z(f)$.
Therefore,
\[
d_m\bigl(Z(f),Z(S(\alpha)f)\bigr) \leq |\alpha| \ \ \ \ \text{and} \
\ \ \ d_m\bigl(Z(\phi_n),Z(S(\alpha)\phi_n)\bigr) = |\alpha|.
\]
Hence $K_m\bigl(S(\alpha)\bigr) = |\alpha|$. The same argument can
be used for $K_F$.

Since $\varpi\bigl(S(\alpha)\bigr) = S(\alpha)\phi_n$, by
Theorem~\ref{tpa2K} we have
\[
K_h\bigl(S(\alpha)\bigr) = \varrho\bigl[ S(\alpha)\phi_n\bigr] =
|\alpha|.
\]
By Proposition~\ref{ps}~\!(\ref{psfn}),
\[
K_H\bigl(S(\alpha)\bigr) = \max\bigl\{K_h
\bigl(S(\alpha)\bigr),K_h\bigl(S(-\alpha)\bigr)\bigr\} = |\alpha|.
\qedhere
\]
\end{proof}

\begin{proposition} \label{estha}
Let $\gamma \in \nC$ and consider $H_k(\gamma) \in \cD(\cP_n)$ for
$k = 1,\ldots,n$. Then
\begin{equation*} 
K_h\bigl(H_k(\gamma)\bigr) = \sqrt[k]{\frac{|\gamma| \,
n!}{(n-k)!}}.
\end{equation*}
\end{proposition}
\begin{proof}
By the definition of $\varpi$,
$\varpi\bigl(H_k(\gamma)\bigr)=H_k(\gamma)\phi_{n}$ and
 \[
\bigl(H_k(\gamma)\phi_{n}\bigr)(z) = \phi_n(z) -\gamma \,
\phi_{n-k}(z) = \frac{z^{n-k}}{n!}\left(z^k -
\frac{\gamma\,n!}{(n-k)!} \right) .
 \]
Consequently,
 $
\varrho\bigl[H_k(\gamma)\phi_{n}\bigr] = \sqrt[k]{|\gamma| \,
n!/(n-k)!},
 $
and the proposition follows from Theorem~\ref{tpa2K}.
\end{proof}

\section{Estimates} \label{ses}

Let $a_0, \ldots, a_n \in \nC$, $a_0\neq0$, and let $T =
T(a_0,\ldots,a_n)$ be the corresponding invertible operator in
$\cF(\cP_n)$. In this section we give estimates for the quantities
$K_h(T)$, $K_H(T)$ and $K_F(T)$ in terms of the coefficients
$a_0,\ldots,a_n$.

\begin{proposition} \label{pHi}
Let $\gamma \in \nC$. Then
\begin{equation} \label{Hng-1}
(n!)^{1/n}|\gamma| \leq K_h\bigl(H_1(\gamma)^{-1}\bigr) \leq n
|\gamma|.
\end{equation}
\end{proposition}
\begin{proof}
To prove the first inequality in \eqref{Hng-1} consider the
polynomial
\[
g = \gamma^{n-1}\, \phi_1 + \gamma^{n-2}\, \phi_2 + \cdots +
\gamma\,\phi_{n-1} + \phi_n.
\]
A straightforward computation gives
 $
H_1(\gamma)g = g - \gamma\, g' = \phi_n - \gamma^n.
 $
Hence, $g$ has a root at $z = 0$, while $H_1(\gamma)g$ has all of
its roots on the circle of radius $(n!)^{1/n}|\gamma|$. If we put $f
= H_1(\gamma)g$, then $g = H_1(\gamma)^{-1}f$ and the previous
observation about the location of the roots implies that
\[
d_h\bigl(Z(f),Z\bigl(H_1(\gamma)^{-1}f \bigr) \bigr) \geq
(n!)^{1/n}|\gamma|.
\]
This implies the first inequality in \eqref{Hng-1}.

To prove the second inequality notice that by \eqref{eqHgi},
\begin{equation*}
\bigl(H_1(\gamma)^{-1}\phi_n\bigr)(z) = \sum_{k=0}^n\,\gamma^{n-k}\,
\phi_k(z) = \gamma^n\,\sum_{k=0}^n\, \phi_k\bigl(z/\gamma\bigr).
\end{equation*}
The roots of $\sum_{k=0}^n \phi_k$ have been researched extensively;
see for example \cite{PV} and the references therein. Here we only
need that the root radius of the polynomial $\sum_{k=0}^n \phi_k$ is
smaller than or equal to $n$.  Consequently, $\varrho
\bigl[H_1(\gamma)^{-1}\phi_n\bigr] \leq n\,|\gamma|$ and the second
inequality in \eqref{Hng-1} follows from Theorem~\ref{tpa2K}.
\end{proof}

\begin{corollary} \label{cKH}
Let $\gamma \in \nC$. Then
\[
\ K_H\bigl(H_1(\gamma)\bigr) = K_H\bigl(H_1(\gamma)^{-1}\bigr) = n
|\gamma|.
\]
\end{corollary}
\begin{proof}
The corollary follows from Propositions~\ref{ps}\!~(\ref{psfn}),
\ref{estha} and \ref{pHi}.
\end{proof}
\begin{proposition} \label{ksh}
Let $\gamma \in \nC$. Then
\[
K_F\bigl(H_1(\gamma)\bigr) \leq n^2 |\gamma|.
\]
\end{proposition}
\begin{proof}
Let $f$ be a non-constant polynomial in $\cP_n$ and $\gamma \in
\nC$. Consider the set
\begin{equation*}
\Omega = \bigcup_{w \in Z(f)} \!
\nD\biggl(\!w+\frac{n\,\gamma}{2},\frac{n\,|\gamma|}{2}\!\biggr).
\end{equation*}
By \cite[Corollary 5.4.1(iii)]{RS} we have
$Z\bigl(H_1(\gamma)f\bigr) \subset \Omega$, and in each connected
component $\Omega_1, \Omega_2, \ldots, \Omega_k$ of $\Omega$ the
polynomials $f$ and $H_1(\gamma)f$ have the same number of zeros,
counted according to their multiplicities. Therefore
\begin{equation*}
d_F\bigl(Z(f),Z(H_1(\gamma)f)\bigr) \leq \max \bigl\{
\diam(\Omega_j) : j=1,\dotsc,k \bigr\}  .
\end{equation*}
Since for each $j = 1,\dotsc,k$,
\begin{equation*}
\diam(\Omega_j) \leq 2n\frac{n\,|\gamma|}{2} = n^2\, |\gamma|,
\end{equation*}
we conclude that
 $
K_F(H_1(\gamma)) \leq n^2 \,|\gamma|.
 $
\end{proof}

\begin{theorem} \label{t13}
Let $a_0,\ldots,a_n\in\nC$, $a_0\neq0$. Let $T=T(a_0,\ldots,a_n)$ be
the corresponding invertible operator in $\cF(\cP_{n})$. Let
$\gamma_1, \ldots, \gamma_n$ be the roots of
\begin{equation*} 
a_0\, z^n + a_1\, z^{n-1}  + \cdots + a_{n-1}\,z + a_n
\end{equation*}
counted according to their multiplicities. Then
\begin{align} \label{t13e1}
 K_h(T) &\leq n \bigl( |\gamma_1| + \cdots + |\gamma_n|   \bigr) \,
, \\ \label{t13e2} %
 K_H(T) & \leq n \bigl( |\gamma_1| + \cdots + |\gamma_n|
\bigr) \, , \\ \label{t13e3}
 K_F(T) &\leq n^2 \bigl( |\gamma_1| + \cdots
                + |\gamma_n| \bigr).
\end{align}
\end{theorem}
\begin{proof}
Lemma~\ref{l3} implies that $\deg(Tf) = \deg(f)$ for all $f \in
\cP_n$. Thus $K_F(T)$ is defined and
$T\bigl(\cP_n\!\setminus\!\cP_0\bigr) = \cP_n\!\setminus\!\cP_0$. By
repeated application of statements (\ref{psa}),(\ref{psen}), and
(\ref{pse}) of Proposition~\ref{ps} to \eqref{ph21} we get
\begin{equation*}
K(T) \leq K\bigl(H_1(\gamma_1)\bigr)+ \cdots +
K\bigl(H_1(\gamma_n)\bigr),
\end{equation*}
where $K$ can be either $K_h$, $K_H$, or $K_F$. Now \eqref{t13e1}
follows from Proposition~\ref{estha}, \eqref{t13e2} follows from
Corollary~\ref{cKH}, and \eqref{t13e3} follows from
Proposition~\ref{ksh}.
\end{proof}

\begin{remark}
The problem of estimating the sum of the absolute values of all the
roots of a given polynomial was studied by Berwald, see
\cite[Theorem~2.3]{MT}.
\end{remark}

\begin{proposition} \label{krks}
Let $\,T$ be an invertible operator in $\cF(\cP_n)$. Then
\begin{equation} \label{krkseq}
K_H(T) \leq K_F(T) \leq (e \, n^3 \ln n)\, K_H(T) .
\end{equation}
\end{proposition}
\begin{proof}
The first inequality was proved in Proposition~\ref{ps}.  Let us now
fix $T \in \cF(\cP_{n})$ and assume that
\[
T = I + \frac{\alpha_0}{1!}\, D + \cdots  +
\frac{\alpha_{n-1}}{(n-1)!}\, D^{n-1} + \frac{\alpha_n}{n!}\, D^n.
\]
By the known estimate \cite[(8.1.12)]{RS}, all the roots $\gamma_j$
of the polynomial
\[
\alpha_0\, z^n + \alpha_1\, z^{n-1} + \cdots +
\frac{\alpha_{n-1}}{(n-1)!}\,z + \frac{\alpha_n}{n!}
\]
are smaller than $\sum_j \bigl(|\alpha_j|/j!\bigr)^{1/j}$. By
\eqref{t13e3} we then have
\begin{equation} \label{krks1}
K_F(T) \leq n^3 \sum_{j=1}^n \,
\left(\frac{|\alpha_j|}{j!}\right)^{1/j}.
\end{equation}

On the other hand, by Theorem~\ref{tpa2K} we have
\begin{equation} \label{krks2}
K_h(T) = d_h\bigl(Z(\phi_{n}), Z(T(\phi_{n})) \bigr) =
\varrho[T\phi_{n}].
\end{equation}
A lower estimate for $\varrho[T\phi_{n}]$, where
\[
\bigl(T\phi_{n}\bigr)(z) = \frac{1}{n!} \biggl( \alpha_n +
\binom{n}{1}\alpha_{n-1} z + \ldots
  + \binom{n}{n-1}\alpha_1 z^{n-1} + z^n\biggr),
\]
is given by another classical inequality (see \cite[(8.1.1)]{RS}):
\begin{equation} \label{krks3}
\max_{1\leq j\leq n} |\alpha_j|^{1/j} \leq \varrho[T\phi_{n}].
\end{equation}
Combining \eqref{krks1}, \eqref{krks3} and \eqref{krks2} and letting
$\mu =\sum_{j=1}^n (j!)^{-1/j}$ gives
\begin{equation*}
K_F(T)\leq n^3 \mu \, \max_{1\leq j\leq n} |\alpha_j|^{1/j}\\
      \leq n^3 \mu \, K_h(T)
      \leq n^3 \mu \, K_H(T).
\end{equation*}

To estimate $\mu$, note that for $j>1$ we have
\[
 \left(\frac{1}{j!}\right)^{1/j} < \frac{e}{j} < \frac{e}{j-\frac{1}{2}} <
 e\ln\frac{j}{j-1},
\]
where the first inequality easily follows from Stirling's
approximation \cite[p.~183]{Mit}, while the last one is a special
case of \cite[3.6.17]{Mit}. Since $1+\frac{1}{\sqrt{2}}<e\ln 2$,
adding up gives
\[
 \mu = \sum_{j=1}^n \left(\frac{1}{j!}\right)^{1/j} < e\ln n
\]
and thus, finally,
 $
 K_F(T) \leq (e\, n^3 \ln n) K_H(T)
 $
as we needed to prove.
\end{proof}

The factor $e\, n^3 \ln n$ in the second inequality of
\eqref{krkseq} is plausibly far from being best possible, but we
have not pursued this line of enquiry.

\section{The main theorem} \label{smt2}

We can now finally give a comprehensive and compact answer to the
question posed in the Introduction.

\begin{theorem} \label{tc2}
Let $\,T \in \cL(\cP_n), \, T\neq 0$. The following statements are
equivalent.
\begin{enumerate}[{\rm (a)}]
\item \label{tc2c}
$T \in \cF(\cP_n)$ and $T$ is invertible.
\item \label{tc2b1}
$T$ is a Grace operator.
\item \label{tc2b}
 $K_h(T) < + \infty$.
\item \label{tc2a}
 $K_m(T) < + \infty$.
\item \label{tc2f}
$T$ is invertible and $TD=D\,T$.
\item \label{tc2d}
 $K_H(T) < + \infty$.
\item  \label{tc3e}
$\deg(Tf) = \deg(f)$ for all $f \in \cP_n$ and $K_F(T) < + \infty$.
\end{enumerate}
\end{theorem}
\begin{proof}
The statements (\ref{tc2c}) through (\ref{tc2f}) are equivalent by
Theorem~\ref{tc1}.

The implication (\ref{tc2c})$\Rightarrow$(\ref{tc3e}) follows from
Proposition~\ref{invert}, equalities \eqref{eqD1} and
Theorem~\ref{t13}. Proposition~\ref{ps}~\!(\ref{psg}) yields
(\ref{tc3e})$\Rightarrow$(\ref{tc2d}) and
(\ref{tc2d})$\Rightarrow$(\ref{tc2b}) follows from
Proposition~\ref{ps}~\!(\ref{psc}). The theorem is proved.
\end{proof}

As $\cF(\cP_n)$ is a commutative subalgebra of $\cL(\cP_n)$ we have
the following corollary.

\begin{corollary}
If $\,T, \, T_1 \in \cL(\cP_n)\!\setminus\!\{0\}$ satisfy any of the
equivalent
conditions~{\rm~(\ref{tc2c})\nolinebreak-\nolinebreak(\ref{tc3e})}
in Theorem~{\rm~\ref{tc2}}, then $T$ and $T_1$ commute.
\end{corollary}

\begin{remark}
By Proposition~\ref{invert} the statement (\ref{tc2c}) in
Theorem~\ref{tc2} is equivalent to $T^{-1}$ is invertible and
$T^{-1} \in \cF(\cP_n)$. Therefore, the statements about the
operator $T$ in Theorem~\ref{tc2} are equivalent to the
corresponding statements about the operator $T^{-1}$.
\end{remark}

\section{Examples} \label{se}

In the next two examples we consider the impact of the operators
from $\cF(\cP_n)$ on the roots of polynomials in $\cP_n$ for $n = 1,
2$, where the computations can be carried out in all detail. While
the situation for $n=1$ is as trivial as expected, the case $n=2$
gives a good insight into why direct calculations for larger
matrices are bound to be unwieldy.

\begin{example}
In ${\cP}_1$ we have
\[
T(a_0,a_1)(\phi_1 - w \phi_0) = a_0 \phi_1 - (a_0 w - a_1)\phi_0.
\]
Provided that $a_0\neq 0$, the polynomial $T(a_0,a_1)(\phi_1 - w
\phi_0)$ has a root at $w-a_1/a_0$. Thus, the operator $T(a_0,a_1)$
shifts all the roots by exactly $a_1/a_0$. This also follows from
$T(a_0,a_1) = a_0 S(a_1/a_0)$; see \eqref{eqSa}.
\end{example}

\begin{example}
Let $w_1,w_2$ be arbitrary complex numbers.  For the operator
$T(a_0,a_1,a_2) \in \cF(\cP_2)$ with $a_0\neq 0$ we have
 \begin{multline*}
\Bigl(T(a_0,a_1,a_2)\bigl(
     (\phi_1 - w_1 \phi_0)(\phi_1 - w_2 \phi_0)\bigr)\!\Bigr)(z) \\
 = a_0 \left( z^2 + 2 \left(
\frac{a_1}{a_0} - \frac{w_1+w_2}{2} \right) z
+w_1w_2-\frac{a_1}{a_0}(w_1+w_2) + \frac{a_2}{a_0} \right)
 \end{multline*}
and so the roots of
 $\,T(a_0,a_1,a_2)\bigl((\phi_1 - w_1 \phi_0)(\phi_1 - w_2 \phi_0)\bigr)$
 are
\begin{equation} \label{roots}
\frac{w_1+w_2}{2} - \frac{a_1}{a_0}\pm\sqrt{\left(
\frac{w_1-w_2}{2}\right)^{\!\!2} + \left( \frac{a_1}{a_0}
\right)^{\!\!2}-\frac{a_2}{a_0}} \, .
\end{equation}
In this case we see that, as predicted by our main result, {\em
both} of the new roots are uniformly close to the original
$w_1,w_2$. Also, we see that if all the constants are real and
$a_1^2\geq a_0a_2$ then an invertible operator in $\cF(\cP_2)$
``sends real roots into real roots''. To check the statement on the
uniform displacement of the roots, do as follows:

Call the roots in \eqref{roots} $z_1$ (with $+$) and $z_2$ (with
$-$). Since we are dealing with complex numbers, let us think of the
square root as having one definite value (then the $\pm$ takes care
of the second root). Also, to simplify the algebra, define
\[
\delta_1 := \frac{a_1}{a_0} \ \ \ \ \ \ \text{and} \ \ \ \ \ \
\delta_2 := \left(\frac{a_1}{a_0} \right)^{\!\!2}-\frac{a_2}{a_0}.
\]
So, for instance, we have
\[
|z_1-w_1| \leq |\delta_1| +
\left|\sqrt{\left(\frac{w_1-w_2}{2}\right)^{\!\!2}+\delta_2} -
\frac{w_1-w_2}{2}\right|
\]
and similarly we can estimate $z_2-w_1$. Now, we have
\[
\left(\!\!\sqrt{\!\left(\frac{w_1-w_2}{2}\right)^{\!\!2}+\delta_2} -
\frac{w_1-w_2}{2}\!\right)\!\!\!\left(\!\!
\sqrt{\!\left(\frac{w_1-w_2}{2}\right)^{\!\!2}+\delta_2} +
\frac{w_1-w_2}{2}\!\right) = \delta_2
\]
and so we see that at least one of the factors on the left hand side
has modulus less than or equal to $\sqrt{|\delta_2|}$. By what noted
above, this means that
\[
\min\bigl\{|z_1-w_1|,|z_2-w_1|\bigr\} \leq |\delta_1| +
\sqrt{|\delta_2|}
\]
which (together with the similar estimates derived with $w_2$
instead of $w_1$) translates into the statement that the roots of
\[
T(a_0,a_1,a_2)\bigl((\phi_1 - w_1\phi_0)(\phi_1 - w_2\phi_0) \bigr)
\]
are not farther away from $w_1$ and $w_2$ than the uniform quantity
\[
\left|\frac{a_1}{a_0}\right| +
\sqrt{\left|\left(\frac{a_1}{a_0}\right)^{\!\!2} -
\frac{a_2}{a_0}\right|}.
\]
\end{example}

\begin{example}
This example shows that $K_H(T) < K_F(T)$. It was found with the
help of Mathematica. We consider $\cP_3$ and the operator $T$ and
its inverse given by
\[
T = I +\frac{2}{3}\,D + \frac{2}{9}\, D^2 -\frac{4}{27}\, D^3, \ \ \
T^{-1} = I -\frac{2}{3}\,D+ \frac{2}{9}\, D^2 +\frac{4}{27}\, D^3.
\]
Set $f(z)=(z-1)^2(z+1)$ and $g(z)=(z+1)^2(z-1)$. Then $Tf = g$ and
$d_F(Z(f),Z(g)) = 2$.  Mathematica calculated that
\begin{equation} \label{eqnum}
\varrho\bigl[T\phi_3\bigr] =
\varrho\bigl[T^{-1}\phi_3\bigr]=\frac{2}{3}\bigl(1+\sqrt[3]{2}\,\bigr)
\approx 1.506614.
\end{equation}
Therefore
\[
K_H(T) = \frac{2}{3}\bigl(1+\sqrt[3]{2}\,\bigr) < 2 \leq K_F(T)
\]
Since both $d_F\bigl(Z(\phi_3),Z(T\phi_3)\bigr)$ and
$d_F\bigl(Z(\phi_3),Z(T^{-1}\phi_3)\bigr)$ are equal to the number
in \eqref{eqnum}, this example also shows $K_F(T)$ is not related to
$\phi_3$ in the sense of Corollary~\ref{cKHe}. We leave it as an
open problem to explore the existence and uniqueness of a polynomial
$f$ such that $K_F(T) = d_F\bigl(Z(f),Z(Tf)\bigr)$.
\end{example}

\begin{example}
Several easy examples of ``bad behavior'' from an operator $T$ which
does not satisfy condition (\ref{tc1c}) in Theorem~\ref{tc1} can be
obtained as follows.

Define $T$ by $T(f):=f + f(0) \phi_n$ for $f\in \cP_n$. Then $T$ is
clearly linear, invertible and and does not commute with $D$. Thus,
by Proposition~\ref{commutant}, $T\not\in\cF(\cP_n)$. If
$w\in\nC\setminus\!\{0,1\}$ is arbitrary, then $Z\bigl(T(\phi_n - w
\phi_0)\bigr)$ consists of all the $n$-th roots of $n! w/(1-w)$,
while all the roots of $\phi_n - w \phi_0$ have modulus
$\sqrt[n]{|w|}$, thus ensuring that the distance (with respect to
each of the distances introduced in Section~\ref{sd}) between the
two multisets of roots can be as large as desired (just let $w$
approach $1$).

The operator $D$ belongs to $\cF(\cP_n)$ but it is not invertible.
With the same polynomials as above we have $Z\bigl(D(\phi_n - w
\phi_0)\bigr) = \{0,\ldots,0\}$ ($n-1$ zeros) again showing that the
distance between $Z\bigl(D(\phi_n - w \phi_0)\bigr)$ and $Z(\phi_n -
w \phi_0)$ (with respect to each of the distances, except $d_F$
which is not defined) can be as large as we wish.

The operator defined by $T(f) = f + (a_{n-1}/n) \phi_0$ for
$f(z)=a_0+a_1z+\cdots+a_nz^n$, is clearly linear, invertible and its
matrix with respect to the basis $\{\phi_0,\ldots,\phi_n\}$ of
$\cP_n$ is upper triangular. This operator does not belong
$\cF(\cP_n)$ since it does not commute with $D$. Polynomials that
show that none of (\ref{tc2a}), (\ref{tc2b}), (\ref{tc2d}) and
(\ref{tc3e}) of Theorem~\ref{tc2} hold are $(\phi_1 - w \phi_0)^n$,
$w \in \nC$. Namely, the roots of $T\bigl((\phi_1 - w
\phi_0)^n\bigr)$ are on the circle centered at $w$ with radius
$\sqrt[n]{|w|}$, and $Z\bigl((\phi_1 - w
\phi_0)^n\bigr)=\{w,\ldots,w\}$ ($n$ times).
\end{example}

\noindent {\sc Acknowledgment}\, The authors thank the referees for
pointing out several related references, the article \cite{S} in
particular.  This resulted in numerous improvements.

\end{document}